\documentclass[journal]{IEEEtran}
\usepackage{cite}
\usepackage{amsmath}
\renewcommand{\vec}[1]{\boldsymbol{#1}}
\usepackage{setspace}
\usepackage{bm}
\usepackage{amssymb}
\usepackage{makecell}
\usepackage{graphicx}
\usepackage{url}
\usepackage{tabularx}
\usepackage{multirow}
\usepackage{float}
\usepackage{slashbox}

\begin{document}

\title{Robust Defense Strategy for Gas-Electric Systems Against Malicious Attacks}

\author{Cheng~Wang, Wei~Wei,~\IEEEmembership{Member,~IEEE}, Jianhui~Wang,~\IEEEmembership{Senior Member,~IEEE}, Feng~Liu,~\IEEEmembership{Member,~IEEE}, Feng~Qiu,~\IEEEmembership{Member,~IEEE}, Carlos~M.~Correa-Posada,~\IEEEmembership{Student~Member,~IEEE}, Shengwei~Mei,~\IEEEmembership{Fellow,~IEEE}
% <-this % stops a space
\thanks{This work is supported by the National Natural Science Foundation of China (51321005). The work of J. Wang and F. Qiu is sponsored by the U. S. Department of Energy. (\textit{Correspond to: Wei Wei})} % <-this % stops a space
\thanks{C. Wang, W. Wei, F. Liu, and S. Mei are with Department of Electrical Engineering, Tsinghua University, 100084 Beijing, China. (c-w12@mails.tsinghua.edu.cn; wei-wei04@mails.tsinghua.edu.cn).}
\thanks{J. Wang and F. Qiu are with Argonne National Laboratory, Argonne, IL 60439, USA (jianhui.wang@anl.gov; fqiu@anl.gov).}
\thanks{C. M. Correa-Posada is with the Colombian System Operator XM, Medellín, Colombia (cmcorrea@xm.com.co).}}

\maketitle

% As a general rule, do not put math, special symbols or citations
% in the abstract or keywords.
\begin{abstract}
—This paper proposes a methodology to identify and protect vulnerable components of connected gas and electric infrastructures from malicious attacks, and to guarantee a resilient operation  by deploying valid corrective actions (while accounting for the interdependency of gas pipeline network and power transmission network). The proposed mathematical formulation reduces to a tri-level optimization problem, where the lower level is a multiperiod economic dispatch of the gas-electric system, the middle level distinguishes the most threatening attack on the coupled physical infrastructures, and the upper level provides optimal preventive decisions to reinforce the vulnerable components and increase the system resilience. By reformulating the lower level problem as a mixed integer linear programming (MILP), a nested column-and-constraint generation (C\&CG) algorithm is developed to solve the min-max-min model. Case studies on two test systems demonstrate the effectiveness and efficiency of the proposed methodology.
\end{abstract}

% Note that keywords are not normally used for peerreview papers.

\section*{Nomenclature}
Most of the symbols and notations used throughout this paper are defined below for quick reference. Others are defined following their first appearances, as needed.

\subsection{Sets and Indices}
\begin{IEEEdescription}[\IEEEusemathlabelsep\IEEEsetlabelwidth{$aaaaa$}]
\item[$t \in T$]     Time periods
\item[$g \in G$]     Traditional units
\item[$n \in N$]     Gas-fired units
\item[$i \in I$]     Power grid nodes
\item[$l \in L$]     Power grid lines
\item[$d \in D$]     Power grid loads
\item[$w \in W$]     Gas wells
\item[$s \in S$]     Gas storage tanks
\item[$c \in C$]     Gas compressors
\item[$j \in J$]     Gas network nodes
\item[$y \in Y$]     Gas network lines
\item[$e \in E$]     Gas loads
\end{IEEEdescription}

\subsection{Parameters}
\begin{IEEEdescription}[\IEEEusemathlabelsep\IEEEsetlabelwidth{$aaaaaaaaaaa$}]
\item[$X$]    Defense budget
\item[$A$]    Attack budget
\item[$P_{g}^{min}/P_{g}^{max}$]  Output range of traditional units
\item[$P_{n}^{min}/P_{n}^{max}$]  Output range of gas-fired units
\item[$R_{+}^g/R_{-}^g$]  Ramp up/down limits of traditional units
\item[$R_{+}^n/R_{-}^n$]  Ramp up/down limits of gas-fired units
\item[$F_l$]  Power transmission line capacity
\item[$p_{dt}$]    Power load demand
\item[$\pi_{l}$]   Reactance of transmission line
\item[$f_{g}(\cdot)$]  Generation cost of traditional units
\item[$Q_d$]  Not-served  load cost coefficient
\item[$u_{gt},u_{nt}$]    Generator status
\item[$\tau^u_j,\tau^l_j$]  Gas pressure range
\item[$q_{s,max}^{in},q_{s,max}^{out}$]  Maximum gas storage in/out range
\item[$r_s^u,r_s^l$]  Maximum/minimum gas storage level
\item[$\beta_n$]   Gas-electric conversion factor
\item[$\gamma_c$]  Compression factor of the compressor
\item[$\phi_y$]    Weymouth equation coefficient
\item[$Q_w,Q_e,Q_s$]   Gas production/not-served/storage cost
\item[$q_{et}$]     Gas load
\item[$Temp$]       Temperature
\item[$Z_{y}$]      Compression factor of the pipeline
\item[$Line_{y}$]   Length of the pipeline
\item[$R_y$]        Diameter
\item[$\rho_0$]     Gas density in standard condition
\item[$\mu$]        Specific gas constant
\item[$F_y$]        Pipeline friction coefficient
\item[$\lambda$]    Unit transformation constant
\end{IEEEdescription}

\subsection{Variables}
\begin{IEEEdescription}[\IEEEusemathlabelsep\IEEEsetlabelwidth{$aaaaaaa$}]
\item[$x_l$]              Defend resources of power transmission line
\item[$x_y,x_n,x_c$]      Defend resources of pipeline/con-line/com-line
\item[$a_l$]              Attack resources of power transmission line
\item[$a_y,a_n,a_c$]      Attack resources of pipeline/con-line/com-line
\item[$z_{gt}/z_{nt}$]    Generator cut-out/connection decision
\item[$p_{gt},p_{nt}$]    Generation output
\item[$\Delta p_{dt}$]    Not-served power load
\item[$\theta$]           Phase angle
\item[$pf_{lt}$]          Power flow
\item[$q_{wt}$]           Gas well output
\item[$q_{st}^{in}/q_{st}^{out}$]   Gas storage in/out rate
\item[$r_{st}$]           Stored gas amount
\item[$q_{yt},q_{ct}$]    Gas flow
\item[$v_{jt}$]           Gas pressure square
\item[$\Delta q_{et}$]    Not-served gas
\end{IEEEdescription}

\section{Introduction}

\IEEEPARstart{O}vER the past decade, technology breakthroughs in horizontal drilling and hydraulic fracturing have fostered the production of natural gas trapped within shale rock formations, and natural gas prices have dropped sharply, creating a growing demand for this energy resource from the electric power industry. However, the growing interdependency between the pipeline network for natural gas and power systems through gas-fired power plants also imposes remarkable challenges for the coordinated planning and operation of both physical systems \cite{Shahidehpour_Gas}, which has been well studied in  \cite{Cong_Dispatch,Xiaping_Planning,OGPF}.

Due to its fast response, gas-fired units are usually managed to compensate, in real-time, the discrepancy between generation and load, especially in the presence of a high penetration of renewable energies. Recently, many interesting and inspiring works have been found. A method for stochastic unit commitment is proposed in \cite{Alabdulwahab_SUC_Gas_Power} for jointly operating natural gas and power systems to mitigate the variability of wind generation. A look-ahead robust scheduling model of a gas-electric system, which takes the natural gas congestion into account, is investigated in \cite{Cong_Look_ahead_congestion}. The importance of modeling the infrastructure for gas delivery in a power system operation with high wind power penetrations is studied in \cite{Devlin_gas_importance}. To absorb excessive wind power, power-to-gas (P2G) technology is comprehensively discussed in \cite{Bailera_P2G} and \cite{Qadrdan_P2G}, which allows a bi-directional interchange of energy, and further enhances the interdependency between natural gas systems and power systems.

Outside the normal operating condition, there could be deliberate sabotage against these critical infrastructures, which might cause catastrophic contingencies such as line tripping and generator outage. The vulnerability assessment for power systems alone is extensively studied in \cite{Salmeron_AD_terrorist, Salmeron_AD_interdiction, Arroyo_AD, Yuan_DAD, Arroyo_Tri_Planning}. The resilience of natural gas networks during conflicts, crises, and disruptions were studied in \cite{Carvalho}. In \cite{Khodayar_gas_AD}, a bilevel linear programming model is proposed to assess the vulnerability of microgrids with multiple energy carriers and then deploy redispatch after deliberate interdictions.  In \cite{Yuan_DAD}, a tri-level programming model is devised to determine the optimal defensive strategy for protecting critical power system components before an attack. A column-and-constraint generation (C$\&$CG) method is suggested to solve the tri-level optimization model. It becomes apparent that deploying preventive defenses could be one of the most effective ways to enhance system reliability.

In fact, from a robust optimization perspective, both nodal power injection uncertainty and unexpected equipment failure can be regarded and modeled as some kind of attack, yielding a min-max-min problem, such as the robust unit commitment problem considered in \cite{Nest_CCG}. This paper proposes a systematic methodology to identify the most threatening attacks against the gas-electric system while taking both preventive and corrective actions into account. The contributions are twofold:

1) A defender-attacker-defender (D-A-D) formulation for protecting critical components of the interconnected gas-electric system. The decision variables of the upper level and lower level are the preventive and corrective actions of the defender, say, the system operator. They strive to minimize the total operation cost in the worst-case disruption, including the penalty of not-served power and gas.  The decision variables of the middle level are controlled by a virtual attacker, who seeks the most disruptive sabotage  against the system. Our formulation is different from existing ones in the following ways. Compared with \cite{Khodayar_gas_AD}, we incorporate the model of reinforcing system components before the attack, which is important for enhancing system resilience. In other words, the formulation in \cite{Khodayar_gas_AD} yields an attacker-defender (A-D) model with a bilevel structure. Compared with \cite{Yuan_DAD}, we model the operation of a natural gas system, which gives rise to a nonlinear and non-convex optimization problem due to the presence of Weymouth equation. Moreover, we consider a multiperiod economic dispatch in the lower level, which captures the subsequent effects of the attack over time; we allow the generators to cut in and out and prohibit over generation after the attack occurs, providing additional flexibility for corrective control. However, this introduces binary variables in the lower level, making the problem much more challenging to solve. This leads to the second contribution.

2) A nested C\&CG algorithm for solving the D-A-D model. The linear min-max-min program in \cite{Yuan_DAD} can be solved by the C\&CG algorithm proposed in \cite{Zeng_CCG}. However, the non-convexity introduced by the gas pressure equation and the discontinuity brought by committing generation units prevents dualizing the lower level problem and complicates the computation. The method in \cite{Khodayar_gas_AD} ignores the discontinuity from committing generators and uses a linear approximation of the Weymouth equation, so as to build a linear lower level problem. In the proposed method, the Weymouth equation is first replaced by its piecewise linear approximation with mixed integer linear form  \cite{Carlos_gas_opf}, yielding a mixed integer linear programming (MILP) lower level, then the nested C\&CG algorithm proposed in \cite{Zeng_CCG_binary} is applied to solve the D-A-D model with discrete recourse decisions. Due to the modeling capacity of MILP, our method is versatile  enough to take into account a broader class of corrective actions in response to the attack.

Although the proposed model currently does not consider the uncertainty of renewable generation, there is no practical difficulty in accounting for it by using a method similar to \cite{Nest_CCG}. Specifically, one can derive the identical compact model to Section II.C if a polyhedral uncertainty set is adopted. The remaining part of the paper is organized as follows. Section II presents the mathematical formulation. Section III derives the solution methodology. Section IV provides numerical results on two test systems to validate the proposed model and algorithm. Finally, Section V presents the conclusion.

\section{Problem Formulation}

\subsection{Defender-Attacker-Defender Formulation}
To evaluate vulnerable components of a networked system, it is natural to introduce a virtual attacker who seeks the most serious disruption to a series of system components, while the defender can take measures before and after the attack that enhance the vulnerable equipment or respond to an attack by redispatching available resources, leading to a tri-level D-A-D model. Before the mathematical formulation is presented, some prerequisite assumptions and simplifications to facilitate model formulation are stated as follows:

1) For modeling the gas system, we assume: i. the gas system operates in steady-state, which means the line pack (and thus the pressure dynamics of gas flow) are ignored. For the dynamic gas system modeling, refer to \cite{Carlos_gas_opf}; ii. the valve of the pipeline stays open unless it is attacked; iii. we adopt the simplified compressor model in \cite{Wolf_BNG}, for detailed modeling of the compressor, please refer to \cite{Cong_Dispatch}.

2) For modeling the commitment of generators, we assume all generators are fast-response units, and that the cut-out/ connection process is instantly completed without delay. This assumption is usually valid for gas-fired units. However, our modeling paradigm has no difficulty in modeling units that are not fast-response, as the minimal on-off time requirements can be formulated as mixed integer linear constraints.

3) For the rules of the defense and attack, we assume: i. one component fails completely if it is not defended and is attacked; ii. the component will be fully functional if it is defended, regardless of being attacked or not.

In the D-A-D model, the defender must find and deploy the optimal protection strategy in the first stage (upper level); then the attacker calculates the optimal attack strategy according to the deployment of defense resources and attacks the system in the second stage (middle level); after that, the defender responds to the attack  to minimize the cost and damage, leading to the following tri-level formulation
\begin{spacing}{0.5}
\begin{equation}
\label{obj_upper}
Obj_{upper}=\mathop{\min}\limits_{\Phi_{up}}F_{OC}
\end{equation}

\begin{equation}
\label{def_budget}
s.t.~~~\sum_lx_l+\sum_yx_y+\sum_nx_n+\sum_cx_c\le X
\end{equation}

\begin{equation}
\label{def_binary}
x_l,x_y,x_n,x_c\in \{0,1\}
\end{equation}

\begin{equation}
\label{obj_middle}
Obj_{middle}=\mathop{\max}\limits_{\Phi_{mid}}F_{OC}
\end{equation}

\begin{equation}
\label{atk_budget}
s.t.~~~\sum_la_l+\sum_ya_y+\sum_na_n+\sum_ca_c\le A
\end{equation}

\begin{equation}
\label{atk_binary}
a_l,a_y,a_n,a_c\in \{0,1\}
\end{equation}

\begin{equation}
\label{obj_lower}
Obj_{lower}=\mathop{\min}\limits_{\Phi_{low}}F_{OC}
\end{equation}

\begin{equation}
\label{replace}
s.t.~~b_{\{\cdot\}}=1-a_{\{\cdot\}}+a_{\{\cdot\}}x_{\{\cdot\}},~~~\{\cdot\}=\{l,n,y,c\}
\end{equation}

\begin{equation}
\label{pgt}
z_{gt}u_{gt}P_{min}^g\le p_{gt} \le z_{gt}u_{gt}P_{max}^g
\end{equation}

\begin{equation}
\label{pnt}
z_{nt}b_nu_{nt}P_{min}^n\le p_{nt} \le z_{nt}b_nu_{nt}P_{max}^n
\end{equation}

\begin{equation}
\label{ramp_up_pg}
p_{g,t+1}-p_{gt}\le z_{gt}u_{gt}R_{+}^g+(1-z_{g,t+1} u_{g,t+1})P_{max}^g
\end{equation}

\begin{equation}
\label{ramp_down_pg}
p_{gt}-p_{g,t+1}\le z_{g,t+1}u_{g,t+1}R_{-}^g+(1-z_{gt} u_{gt})P_{max}^g
\end{equation}

\begin{equation}
\label{ramp_up_pn}
p_{n,t+1}-p_{nt}\le z_{nt}u_{nt}R_{+}^n+(1-z_{n,t+1} u_{n,t+1})P_{max}^n
\end{equation}

\begin{equation}
\label{ramp_down_pn}
p_{nt}-p_{n,t+1}\le z_{n,t+1}u_{n,t+1}R_{-}^n+(1-z_{nt} u_{nt})P_{max}^n
\end{equation}

\begin{equation}
\label{load}
0 \le \Delta p_{dt} \le p_{dt}
\end{equation}

\begin{equation}
\label{theta}
-\pi \le \theta_{it} \le \pi
\end{equation}

\begin{equation}
\begin{split}
\label{power_balance}
& \sum_{\{\cdot\}\in \phi_{\{\cdot\}}(i)}p_{\{\cdot\}t}+\sum_{l\in \phi_{O_2}(i)}pf_{lt}-\sum_{l\in\phi_{O_1}(i)}pf_{lt}-\cdots \\ &-\sum_{d\in \phi_d(i)}(p_{dt}-\Delta p_{dt})=0,~~\{\cdot\}=\{g,n\}
\end{split}
\end{equation}

\begin{equation}
\label{power_flow}
-F_{l}\le pf_{lt} \le F_{l}
\end{equation}

\begin{equation}
\label{theta_flow}
\pi_lpf_{lt}=b_{l}(\theta_{i_{1}t}-\theta_{i_{2}t}),~~i_1\in {O_1}(l),i_2\in {O_2}(l)
\end{equation}

\begin{equation}
\label{gas_well}
q_w^l \le q_{wt} \le q_w^u
\end{equation}

\begin{equation}
\label{storage_amount}
r_s^l\le r_{st}=r_{s,t-1}+q_{st}^{in}-q_{st}^{out}\le r_s^u
\end{equation}

\begin{equation}
\label{storage_in}
0 \le q_{st}^{in} \le q_{s,max}^{in}
\end{equation}

\begin{equation}
\label{storage_out}
0 \le q_{st}^{out} \le q_{s,max}^{out}
\end{equation}

\begin{equation}
\label{node_pressure}
(\tau_j^l)^2\le v_{jt} \le (\tau_j^u)^2
\end{equation}

\begin{equation}
\begin{split}
\label{gas_balance}
&\sum_{s\in_{\phi_s(j)}}(q_{st}^{out}-q_{st}^{in})-\sum_{e\in_{\phi_e(j)}}(q_{et}-\Delta q_{et})+ \\
&+\sum_{w\in_{\phi_w(j)}}q_{wt}-\sum_{n\in_{\phi_n(j)}}p_{nt}/\beta_n=\sum_{\{\cdot\}\in \phi_{\{\cdot\}_{1}}(j)}q_{\{\cdot\}t}-\\ &\sum_{\{\cdot\}\in \phi_{\{\cdot\}_{2}}(j)}q_{\{\cdot\}t},~~\{\cdot\}=\{c,y\}
\end{split}
\end{equation}

\begin{equation}
\label{gas_load}
0\le \Delta q_{et}\le q_{et}
\end{equation}

\begin{equation}
\label{weymouth}
q_{yt}|q_{yt}|=\phi_{y}b_y(v_{j_{1}t}-v_{j_{2}t})
\end{equation}

\begin{equation}
\label{wey_co}
\phi_y=\frac{\pi^2\lambda^2R_y^5}{16Line_yF_y\mu TempZ_{y}\rho_0^2}
\end{equation}

\begin{equation}
\label{compressor_pressure}
v_{j_2}\le \gamma_c^2b_cv_{j_1}+(1-b_c)(\tau_{j_2}^u)^2,~~j_1,j_2\in \phi(c)
\end{equation}

\begin{equation}
\label{compressor_flow}
0\le q_{ct} \le bigM_{gf}b_c
\end{equation}

\begin{equation}
\begin{split}
\label{Operation_Cost}
F_{OC}=\sum_{t} &(\sum_gf_g(p_{gt})+\sum_wQ_wq_{wt}+\sum_sQ_sq_{st}+\\ &\sum_dQ_d\Delta p_{dt}+\sum_eQ_e\Delta q_{et})
\end{split}
\end{equation}

\begin{equation}
\label{variable_set_upper}
\Phi_{up}=\{x_l,x_n,x_y,x_c\}
\end{equation}

\begin{equation}
\label{variable_set_mid}
\Phi_{mid}=\{a_l,a_n,a_y,a_c\}
\end{equation}

\begin{equation}
\begin{split}
\label{variable_set_lower}
\Phi_{low}=& \{z_{gt},z_{nt},p_{gt},p_{nt},\theta_{it},pf_{lt},p_{mt},\Delta p_{dt}, \\ & q_{wt},q_{st}^{in},q_{st}^{out},r_{st},q_{yt},q_{ct},v_{jt},\Delta q_{et}\}
\end{split}
\end{equation}
\end{spacing}
\begin{spacing}{2.5}
\end{spacing}
\begin{spacing}{0.95}

In this formulation, the upper level defense problem consists of (\ref{obj_upper})-(\ref{def_binary}). The objective function (\ref{obj_upper}) is to minimize $F_{OC}$. (\ref{def_budget}) is the defense budget constraint where $X$ is a positive integer and (\ref{def_binary}) restricts the decision variables to binary. The middle level attack problem consists of (\ref{obj_middle})-(\ref{atk_binary}). The objective function (\ref{obj_middle}) is to maximize $F_{OC}$. Similarly, (\ref{atk_budget}) is the constraint of the attack budget where $A$ is a positive integer and (\ref{atk_binary}) imposes binary restriction on the attack variables. The lower level problem consists of (\ref{obj_lower})-(\ref{compressor_flow}). The objective function (\ref{obj_lower}) is to minimize $F_{OC}$. (\ref{replace}) represents the operating availability of components after attack, one for normal and zero for failure. (\ref{pgt})-(\ref{compressor_pressure}) constitute the operating constraints of the gas-electric system. For the electricity network, (\ref{pgt}) and (\ref{pnt}) enforce the generation capacity of traditional units and gas-fired units, respectively. (\ref{ramp_up_pg})-(\ref{ramp_down_pn}) are the ramping rate limits of traditional units and gas-fired units, respectively. (\ref{load}) sets the boundary of not-served load. (\ref{theta}) describes the upper and lower phase angle limits of the power grid. (\ref{power_balance}) depicts the power balancing condition. (\ref{power_flow}) is the flow limit for network power. (\ref{theta_flow}) describes the relationship between power flow and phase angle. For gas systems, (\ref{gas_well}) limits the production capacity of gas wells. (\ref{storage_amount}) depicts the capacity for gas storage. (\ref{storage_in}) and (\ref{storage_out}) impose the range of in and out rate of gas storage, respectively. (\ref{node_pressure}) is the pressure limit of each node. (\ref{gas_balance}) is the gas balance requirement. (\ref{gas_load}) declares the boundary of the not-served gas load. (\ref{weymouth}) is the Weymouth equation and it characterizes the relationship between gas flow in a passive pipeline and node pressure. (\ref{wey_co}) defines the coefficient $\beta_y$ in the Weymouth equation. (\ref{compressor_pressure}) implies the pressure relationship between the initial and terminal node of an active pipeline. (\ref{compressor_flow}) states the boundary of gas flow in the active pipeline, in which $bigM_{gf}$ is a large enough positive number. In practice, $bigM_{gf}$ can be chosen as the upper bound of the production capacity of the largest gas well. Particularly, the pipelines with and without compressors are referred to as active and passive pipelines, respectively. (\ref{Operation_Cost}) gives the expression of $F_{OC}$, in which the first two terms represent the production cost of power and gas, the third term represents the cost of storing gas, the last two terms represent the penalty of not-served power and gas. In (\ref{Operation_Cost}), $f_g(\cdot)$ might be quadratic and can be further linearized by using the piecewise linear (PWL) approximation method \cite{Carrion_PWL}. (\ref{variable_set_upper})-(\ref{variable_set_lower}) define the sets of decision variables for upper, middle and lower level problems, respectively.

Compared with the existing models in \cite{Arroyo_AD} and \cite{Yuan_DAD} which minimize the sum of not-served power load in a single period, our work reveals the vulnerability information from both economic and security perspectives. To see this, consider the situation that the system is highly flexible or the attack budget is limited, such that there is no load shedding at all. Then no vulnerable components are detected by the existing methods. However, such a flexible response may involve dispatching costly resources in real time. In this regard, we incorporate the production cost of the gas and power load, the operating cost of storage, and the cost of load shedding into the objective function, which allows the operator to protect vulnerable facilities before an attack occurs in a more economical manner.

Another important distinction is considering over-generation of units after the attack, which is neglected in the literature \cite{Arroyo_AD}, \cite{Yuan_DAD}, \cite{Khodayar_gas_AD}. According to the attacker, a direct attack is to destroy the lines that connect the load with the network. However, an indirect attack, which destroys some lines to cause operation infeasibility even if the generators' output reaches the minimum value in the recourse stage (that is, the recourse action leads generators to an infeasible operating point $P^{min}_g/2$), may cause more catastrophic consequences as some generators have to be cut out, leading to a blackout of the system. From the analysis above, we can see, the unmodeled feature will cause conflict between generator capacity constraint (\ref{pgt})-(\ref{pnt}) and nodal power balance constraint (\ref{power_balance}), leading to a suboptimal attack strategy. Physically, over-generation  may trigger a self-protection module and cause a generator outage, leading to an even more catastrophic consequence than the attack itself. To address this issue, the commitment strategy of generators is considered and is controlled by binary variables $z_{gt}, z_{nt}$ in (\ref{pgt})-(\ref{ramp_down_pn}). Its significance is demonstrated in Section IV.C.
\end{spacing}

\vspace{-9pt}
\subsection{Linearizing the Weymouth Equation}

Clearly, each term in the proposed tri-level formulation is linear except for the gas flow square with sign function in the Weymouth equation (\ref{weymouth}), which can be linearized by PWL approximation as follows \cite{Carlos_gas_opf}.

\begin{spacing}{0.5}
\begin{equation}
\label{gas_flow_1}
\psi(q_{yt})\approx\psi(\Delta_{y1})+\sum_{k\in K}(\psi(\Delta_{y,k+1})-\psi(\Delta_{yk}))\varrho_{ytk}
\end{equation}

\begin{equation}
\label{gas_flow_2}
q_{yt}=\Delta_{y1}+\sum_{k\in K}(\Delta_{y,k+1}-\Delta_{y,k})\varrho_{ytk}
\end{equation}

\begin{equation}
\label{gas_flow_3}
\varrho_{yt,k+1}\le\zeta_{ytk}\le\varrho_{ytk},~\forall k \in K-1
\end{equation}

\begin{equation}
\label{gas_flow_4}
0\le\varrho_{ytk}\le1,~\forall k\in K
\end{equation}
\end{spacing}
\begin{spacing}{1.5}
\end{spacing}
\noindent Where $\Delta_{yk}$ is the piecewise segment for each passive pipeline and $k$ is the corresponding index; $\varrho_{ytk}$ and $\zeta_{ytk}$ are the auxiliary continuous and binary variable, respectively; $\psi$ depicts the relationship between the nonlinear term and gas flow. You can control the error gap of PWL approximation (\ref{gas_flow_1})-(\ref{gas_flow_4}) by choosing the size of $K$ \cite{Carlos_linear_gas}.

\vspace{-9pt}
\subsection{The Compact Model}

For ease of analysis, the tri-level D-A-D model in previous subsections can be written in a compact form as follows:
\begin{spacing}{0.5}
\begin{equation}
\label{obj_upper_mat}
Obj_{upper}=\mathop{\min_{\vec{x}}}\vec{h}^T\vec{y}~~~~~~~~~~~~~~~~~~~~~~~
\end{equation}

\begin{equation}
\label{upper_con_mat}
s.t.~~~\vec{Ax}\le\vec{b}~~~~~~~~~~~~~~~~~~~~
\end{equation}

\begin{equation}
\label{obj_mid_mat}
Obj_{middle}=\mathop{\max_{\vec{a}}}\vec{h}^T\vec{y}~~~~~~~~~~~~
\end{equation}

\begin{equation}
\label{mid_con_mat}
s.t.~~~\vec{Ca} \le \vec{d}~~~~~~~~~~~~
\end{equation}

\begin{equation}
\label{obj_lower_mat}
Obj_{lower}=\mathop{\min_{\vec{y, z}}}\vec{h}^T\vec{y}~~~~
\end{equation}

\begin{equation}
\label{lower_con_mat}
~~~~~~~~~~~s.t.~~~\vec{E(x, a)y+F(x, a)z}\le\vec{g}
\end{equation}
\end{spacing}
\begin{spacing}{2.5}
\end{spacing}

\noindent Where $\vec{x, a, z}$ are binary variable vectors, $\vec{y}$ is a continuous variable vector, and $\vec{A, C, b, d, g}$, and $\vec{h}$ are constant coefficient matrices which can be derived from (\ref{def_budget}), (\ref{atk_budget}), and (\ref{pgt})-(\ref{compressor_flow}). Specifically, $\vec{E(x, a)}$ and $\vec{F(x, a)}$ are variable coefficient matrices and can be derived from (\ref{pgt})-(\ref{compressor_flow}). It is worth mentioning that in its current formulation, we consider only the redispatch cost $\vec {h^T y}$ in the objective function, however, there is no difficulty if one prefers to use a more general cost function in the form of $\vec {h^T_1 x + h^T_2 z + h^T_3 y}$.

\section{Solution Methodology}

In general, several decomposition methods are available to solve a linear min-max-min model such as Benders decomposition \cite{Geoffrion_Benders,RUC-Benders} and column and constraint generation ($C\&CG$) \cite{Zeng_CCG}. However, the presence of binary variables prevents directly dualizing the lower level problem. We adopt a nested C\&CG method in \cite{Zeng_CCG_binary} to solve the D-A-D problem (\ref{obj_upper_mat})-(\ref{lower_con_mat}). The application of the nested C\&CG method has been found in the robust unit commitment with quick-start units \cite{Nest_CCG}, which shows its potential  for the proposed formulation.

\vspace{-6pt}
\subsection{Inner $C\&CG$ for Middle Lower-Level Problem}

With fixed defense strategy $\vec{x^*}$ and lower-level binary variable $\vec{z^*}$, (\ref{lower_con_mat}) yields the following inequality

\begin{spacing}{0.5}
\begin{equation}
\label{reformed_lower_con}
\vec{E(x^*,a)y}\le\vec{g-F(x^*,a)z^*}
\end{equation}
\end{spacing}

\begin{spacing}{1.5}
\end{spacing}
\begin{spacing}{0.95}
Then the linear program consisting of (\ref{obj_lower_mat}) and (\ref{reformed_lower_con}) with $\vec{z=z^*}$ can be replaced by its dual form as follows
\end{spacing}

\begin{spacing}{0.5}
\begin{equation}
\label{reformed_lower_obj_mat}
\mathop{\max_{\vec{\mu}}}(\vec{g-F(x^*,a)z^*})^T\vec{\mu}
\end{equation}

\begin{equation}
\label{reformed_lower_dual_con1}
s.t.~ \vec{E(x^*,a)}^T\vec{\mu}=\vec{h}, \vec{\mu}\le\vec{0}
\end{equation}
\end{spacing}

\begin{spacing}{1.5}
\end{spacing}

\noindent Where $\vec{\mu}$ is the dual variable. The middle lower-level problem (\ref{obj_mid_mat})-(\ref{lower_con_mat}) can be solved by the $C\&CG$ algorithm described below, denoted as $Algo_{inner}$.

\subsubsection*{Step 1} Select an arbitrary feasible attack strategy $\vec{a^*}$, and solve the following problem:
\begin{equation}
\label{CCG_inner_subp}
\begin{gathered}
\mathop{\min_{\vec{y, z}}}\vec{h}^T\vec{y}  \\
s.t.~\vec{E(x^*, a^*)y+F(x^*, a^*)z}\le\vec{g}
\end{gathered}
\end{equation}
Denote the optimal solution by $\vec{y^*,z^*}$, set $LB=\vec{h^T y^*}$, $UB=+\infty$, $o=1$, $\vec{z^{1*}=z^*}$, $O=\{ 1 \}$, $Z=\{\vec{z^{1*}} \}$.
\subsubsection*{Step 2} Solve the following problem

\begin{spacing}{0.5}
\begin{equation}
\label{CCG_inner_obj}
Obj_{middle}=\mathop{\max_{\theta,\vec{a,\mu}}}\theta ~~~~~~~
\end{equation}

\begin{equation}
\label{CCG_inner_con_mid1}
s.t.~\theta\le\vec{(g-F(x^*,a)z^{r*})^T\mu^r},~\forall r\in O, \vec{z^{r*}} \in Z
\end{equation}

\begin{equation}
\label{CCG_inner_con_mid2}
\vec{Ca} \le \vec{d},\ \vec{E(x^*,a)^T\mu^r}\le\vec{h^T},\vec{\mu^r}\le\vec{0},~\forall r \in O
\end{equation}
\end{spacing}

\begin{spacing}{2.0}
\end{spacing}
\noindent Extract the optimal solution $\vec{a^*}$ and $\theta^*$, and update $UB = \theta^*$.

\subsubsection*{Step 3} Solve (\ref{CCG_inner_subp}) with $\vec{a^*}$, obtain the optimal solution $(\vec{z^*,y^*})$ and update $LB={\max}\{LB,\vec{h^T y^*}\}$.

\subsubsection*{Step 4}  If $UB-LB\le\epsilon$, terminate and, return the attack strategy $\vec{a^*}$ and the optimal value $Obj_{middle}^*$. Otherwise, update $o=o+1$,  $O=O\cup o+1$, $\vec{z^{o*}} = \vec{z^*}$, $Z = Z \cup \vec{z^{o*}}$, create new variables $\vec{\mu^{o}}$ with the following constraints,

\begin{spacing}{0.5}
\begin{equation}
\label{CCG_inner_con_CCG1}
\theta\le\vec{(g-F(x^*,a)z^{o*})^T\mu^{o}}
\end{equation}

\begin{equation}
\label{CCG_inner_con_CCG2}
\vec{E(x^*,a)^T\mu^{o}}\le\vec{h^T},\vec{\mu^{o}} \le \vec{0}
\end{equation}
\end{spacing}

\begin{spacing}{1.5}
\end{spacing}
\noindent and go to Step 2.

\begin{spacing}{0.95}
Note that there are bilinear terms in (\ref{CCG_inner_con_mid1}) and (\ref{CCG_inner_con_mid2}) which can be further linearized, and Section III.C discusses them. Then problem (\ref{CCG_inner_obj})-(\ref{CCG_inner_con_mid2}) can be reformulated as a MILP, which is readily solvable by commercial software. By calling $Algo_{inner}$, the optimal attack strategy $a^*$ under fixed defense strategy $x^*$ can be obtained.
\end{spacing}

\vspace{-6pt}
\subsection{Outer $C\&CG$ for Upper Level Problem}

The outer $C\&CG$ algorithm for the upper level problem, denoted as $Algo_{outer}$, identifies the optimal defense strategy under all possible attacks, which proceeds as follows.

\subsubsection*{Step 1} Set $LB=-\infty$, $UB=+\infty$, $\vec{a^{1*}}=0$, $w=1$, $W =\{w\}$, and select a convergence tolerant $\epsilon_{outer}$.
\subsubsection*{Step 2} Solve the following problem:
\begin{gather}
Obj_{upper}=\mathop{\min_{\vec{x,y,z}}}\varphi   \label{CCG_outer_obj} \\
s.t.~\vec{Ax}\le\vec{b},\ \varphi\le\vec{h^Ty^f},~\forall f\in W  \label{CCG_outer_con_1}\\
\vec{E(x,a^{f*})y^f+F(x,a^{f*})z^f}\le\vec{g},~\forall f\in W  \label{CCG_outer_con_2}
\end{gather}
Let the optimal solution be $\vec{x^*}$, the optimal value $\varphi^*$, and update $LB=\mathop{\min}(UB, \varphi^*)$.
\subsubsection*{Step 3} Call $Algo_{inner}$, obtain the optimal solution $\vec{a^*}$ and  optimal value $Obj_{middle}^*$, update $UB=\mathop{\min}(UB,Obj_{middle}^*)$.
\subsubsection*{Step 4} If $UB-LB\le\epsilon_{outer}$, terminate and return the optimal solution $\vec{x^*,a^*,y^*,z^*}$. Otherwise, update $w=w+1$, $W=W\cup w$, $\vec{a^{w*}}=\vec{a^*}$, create new variables $(\vec{y^{w},z^{w}})$ with the following constraints, and go to Step 2.

\begin{spacing}{0.5}
\begin{equation}
\label{CCG_outer_con_1_CCG}
\varphi\le\vec{h^Ty^{w}}
\end{equation}

\begin{equation}
\label{CCG_outer_con_2_CCG}
\vec{E(x,a^{w*})y^{w}+F(x,a^{w*})z^{w}}\le\vec{g}
\end{equation}
\end{spacing}

\begin{spacing}{1.5}
\end{spacing}

Similarly, there are bilinear terms in (\ref{CCG_outer_con_2}) which can be further linearized, and the next subsection discusses them. Then the problem in Step 2 can be transformed into a MILP.

\vspace{-8pt}
\subsection{Linearization of Bilinear Terms}

The bilinear terms in $Algo_{inner}$ and $Algo_{outer}$ include binary production and binary-continuous production. The corresponding linearization methods are shown below.

\subsubsection{Binary production}

The binary production terms exist in (\ref{pnt}), namely $\vec{F(x,a^*)z^r}$ of (\ref{CCG_outer_con_2}). Take $\delta = \kappa \lambda$ for example, where $\kappa$ and $\lambda$ are binary, and their multiplication $\delta$ is another binary variable. It is easy to verify that above equality is equivalent to the following linear constraints
\begin{gather}
\kappa+\lambda-\delta-1\le 0  \label{binary_1}  \\
\delta-\kappa\le 0   \label{binary_2}   \\
\delta-\lambda\le 0   \label{binary_3}
\end{gather}

\subsubsection{Binary-continuous production}

The binary-continuous production terms exist in (\ref{theta_flow}), (\ref{weymouth}), (\ref{compressor_pressure}), namely $\vec{E(x,a^*)y^r}$ of (\ref{CCG_outer_con_2}). Again, take $\delta=\kappa\lambda$ for example, where $\kappa$ is binary and $\lambda$ is continuous. Their multiplication $\delta$ is a continuous variable. It is easy to verify that above equality is equivalent to the following linear constraints
\begin{gather}
\lambda^{min}\kappa\le\delta\le\lambda^{max}\kappa   \label{half_binary_1}\\
\lambda^{min}(1-\kappa)\le\lambda-\delta\le\lambda^{max}(1-\kappa)  \label{half_binary_2}
\end{gather}
where $\lambda^{min}$ and $\lambda^{max}$ are the lower bound and upper bound of $\lambda$, respectively. If they are not known precisely, we can use a large constant $BigM_{con}$ instead. In practice, choosing the value of $BigM_{con}$ can be tricky: if $BigM_{con}$ is too large, the continuous relaxation of the MILP will be rather weak, this will cause unnecessary computational burden in the branch and bound procedure; on the other hand, if $BigM_{con}$ is too small, there will be a duality gap between problem (\ref{reformed_lower_obj_mat})-(\ref{reformed_lower_dual_con1}) and its primal problem, and the reformulation will be incorrect. Due to the physical operating constraints, the bounds of primal variables (say the gas pressure and the power transmission line capacity) are easy to retrieve, and the bounds of dual variables are not instantly apparent.

To reduce the proper value of $BigM_{con}$, the objective function (\ref{Operation_Cost}) is divided by an auxiliary parameter $BigM_{obj}$. A heuristic choice of $BigM_{obj}$ may be the objective value when all loads are not-served. To see this, consider the following primal-dual pair of LPs:
\begin{equation*}
\begin{gathered}
\mbox{P:  } \min ~ \{  c^T x ~|~ Ax \le b \}  \\
\mbox{D: } \max ~\{ b^T y ~| ~ A^Ty = c, y<=0 \}
\end{gathered}
\end{equation*}
Due to physical operation constraints, the feasible region of problem P is a bounded polytope, and a finite optimal value $c^T x^*$ can be found at one of its extreme points $x^*$. According to the duality theory of LP, problem D must have a finite optimal value $b^T y^* = c^T x^*$ with bounded $y^*$, although the feasible region of D may be unbounded. In transformation (\ref{half_binary_1}) and (\ref{half_binary_2}), the $BigM_{con}$ parameter of the dual variable should satisfy  $-BigM_{con} \le y^*$ such that the optimal solution is not influenced. The difficulty lies in the fact that $y^*$ is not known in advance, so we have to use a large enough value, which might be overly pessimistic. However, if the objective of P becomes $c^T x/BigM_{obj}$, the optimal solution $x^*$ of problem P remains the same, and it is easy to verify that $y = y^*/BigM_{obj}$ is an optimal solution of problem D. This indicates that the bound parameter $BigM_{con}$ imposed on the dual variables can be reduced. Section IV.D shows the numeric results.

\section{Illustrative Example}

In this section, we present numerical experiments on two test systems to show the effectiveness of the proposed model and algorithm. The experiments are performed on a laptop with Intel(R) Core(TM) 2 Duo 2.2 GHz CPU and 4 GB memory. The proposed algorithms are implemented on MATLAB with YALMIP toolbox. MILP is solved by Gurobi 6.5, the optimality gap is set as $0.1\%$ without particular mention.

\subsection{6-Bus Power System with 7-Node Gas System}

\begin{figure}
\centering
  \includegraphics[width=0.35\textwidth]{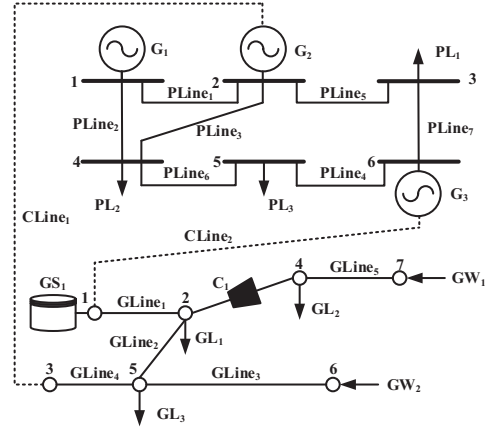}
  \caption{Topology of the hybrid test system.}
  \label{fig:Topology_P6G7}
\end{figure}

\begin{figure}
\centering
  \includegraphics[width=0.35\textwidth]{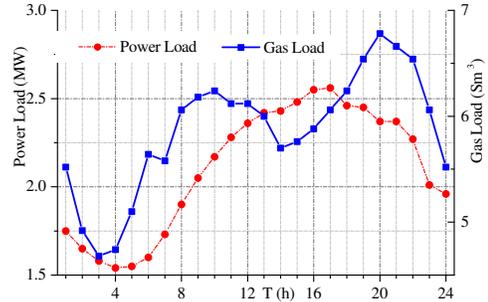}
  \caption{Power and gas demand.}
  \label{fig:gasload}
\end{figure}

Fig. \ref{fig:Topology_P6G7} depicts the topology of the connected infrastructure. It has 2 gas-fired units, 1 traditional unit, 2 gas wells, 1 compressor, 1 gas storage tank, 3 power loads and 3 gas loads. The targets of the attacker are the 7 power transmission lines, 5 passive pipelines, 1 active pipeline and 2 connection lines. Here the connection lines represent the gas pipelines that transport gas to gas-fired units. The parameters of the system and the scheduled unit commitment can be found in \cite{Power6Gas7}. In the following cases, we consider the problem with $T=4$ periods. The power and gas demand profiles are shown in Fig. \ref{fig:gasload}, in which four typical time slots are selected. Details of these time slots are shown in Table \ref{Tab:time plots}. We set $K=8$ in (\ref{gas_flow_1})-(\ref{gas_flow_4}), which means that an eight-segment linear approximation replaces the nonlinear Weymouth equation.
\begin{spacing}{0.1}
\end{spacing}
\begin{table}[h]
\footnotesize
  \centering
  \caption{Details of the Selected Time Slots}\label{Tab:time plots}
  \begin{tabular}{c|c|c|c}
  \hline
  % after \\: \hline or \cline{col1-col2} \cline{col3-col4} ...
  Case & Periods & Power & Gas \\
  \hline
  1 & 2-5 & low & low \\
  \hline
  2 & 9-12 & mid & mid \\
  \hline
  3 & 14-17 & high & mid \\
  \hline
  4 & 19-22 & mid & high \\
  \hline
  \end{tabular}
\end{table}

\subsection{Computational Results}

\begin{figure}
\centering
  \includegraphics[width=0.40\textwidth]{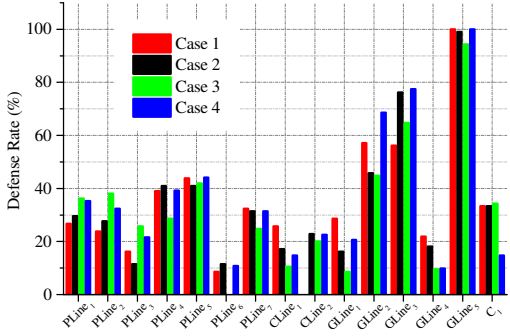}
  \caption{Defense rate of components under 4 cases.}
  \label{fig:Defend_Budget}
\end{figure}

In this subsection, $X=3$ and $A=3$ are selected as the benchmark and the corresponding results are summarized in Table \ref{Tab:case4benchmark}, where \emph{Def}, \emph{Atk}, \emph{OC}, \emph{NP}, \emph{NG} are short for \emph{Defend Resource}, \emph{Attack Resource}, \emph{Operational Cost}, \emph{Not-served Power} and \emph{Not-served Gas}, respectively. As the gas system is radial and the power system is meshed in topology, the production resources of the gas system are more centralized than the power system and more vulnerable to attacks. Thus the defense resources are spent mainly on the gas system. Similar circumstance occurs in the attacker's strategy. The attack priority is the gas system unless the gas load is in low level, that is, case 1. The average power and gas loss rates of the four time plots under the optimal attack strategy are 48.6$\%$ and 6.55$\%$, respectively, assuming the optimal deployment of defense resources. This difference reveals that the failure of the gas system will have a larger impact on the coupled physical system, as the energy from gas system to power system is a one-way flow. Also, we consider both the defense and attack budget variables, whose range varies in the intervals $[1, N+L+C+Y-1]$ and $[1, N+L+C+Y-A]$, respectively. Various tests are carried out under different defense budgets, attack budgets, and time slots. Fig. \ref{fig:Defend_Budget} shows the results. Here we define the defense rate of a component as the number of instances in which it is protected divided by the total number of defense and attack budget combinations, for example, $GLine_{5}$ is defended in 99 out of 105 defense and attack budget combinations in case 3, thus its defense rate is $94.26\%$, as shown in Fig. \ref{fig:Defend_Budget}. The importance of each component is seen by the relative value of the defense rate in Fig. \ref{fig:Defend_Budget}, from which we can observe $GLine_{5}$ is always the most important component, as its defense rate is the highest in all cases. Also, the defense rate of a specific component varies from case to case, which reveals that the defense strategy depends on the load level.

\begin{figure}
\centering
  \includegraphics[width=0.33\textwidth]{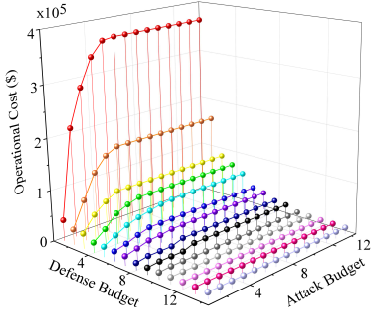}
  \caption{Operational cost under different budgets.}
  \label{fig:Loss}
\end{figure}

To demonstrate the benefit of prior protection, the operation cost under different combinations of defense and attack budget are shown in Fig. \ref{fig:Loss}. from which we can see, with defense budget increasing, the average operation cost under different attack budgets decreases. However, the decrement of operational cost also decreases. Specifically, when the defense budget is 8 or larger, the attack budget has little impact on the operation cost. Moreover, if the attack budget is limited, the operation cost can be well-controlled through a limited defense budget, say $X=4$ in this case, which may provide guidance for defense budgeting.

\begin{table}[h]
\footnotesize
  \centering
  \newcommand{\tabincell}[2]{\begin{tabular}{@{}#1@{}}#2\end{tabular}}
  \caption{Computational Results under the Benchmark}\label{Tab:case4benchmark}
  \begin{tabular}{m{0.5cm}<{\centering}|m{1cm}<{\centering}|m{1cm}<{\centering}|m{1.4cm}<{\centering}|m{0.8cm}<{\centering}|m{0.8cm}<{\centering}}
    \hline
    Case & Def & Atk & OC(\$) & \tabincell{c}{NP\\($MWh$)}& \tabincell{c}{NG\\($Sm^3$)}\\
    \hline
    1 &$GLine_{2}$ $GLine_{3}$ $GLine_{5}$&$GLine_{1}$ $GLine_{4}$ $PLine_{5}$&$6.713\times 10^5$&$664$&$0$\\
    \hline
    2 &$C_{1}$ $GLine_{3}$ $GLine_{5}$&$PLine_{2}$ $PLine_{3}$ $PLine_{4}$&$7.529\times 10^5$&$310$&$2261$\\
    \hline
    3 &$C_{1}$ $GLine_{5}$ $PLine_{5}$&$CLine_{1}$ $PLine_{1}$ $PLine_{2}$&$8.615\times 10^5$&$0$&$3900$\\
    \hline
    4 &$C_{1}$ $GLine_{3}$ $GLine_{5}$&$PLine_{2}$ $PLine_{3}$ $PLine_{4}$&$8.047\times 10^5$&$795$&$0$\\
    \hline
  \end{tabular}
\end{table}

\subsection{Significance of Considering Over-generation}

This subsection discussed the significance of incorporating over-generation. To make the result clear, we assume only power transmission lines can be attacked in the tests. We set $X=2$ and $A=2$, and list the corresponding results in Table~\ref{Tab:no_over_generation}. These results shown that, the sum of the attack strategy is smaller than the attack budget in cases 1-4 if over-generation is not considered, which indicates that any additional attack will cause over-generation and will not be allowed. Further, the optimal defense strategy with or without over-generation attack (OGA) is also quite different, which confirms the need to consider over-generation when deploying corrective actions.

\begin{table}[h]
\footnotesize
  \centering
  \newcommand{\tabincell}[2]{\begin{tabular}{@{}#1@{}}#2\end{tabular}}
  \caption{Computational Result without Over-generation Attack}\label{Tab:no_over_generation}
  \begin{tabular}{m{0.5cm}<{\centering}|m{0.5cm}<{\centering}|m{1cm}<{\centering}|m{1cm}<{\centering}|m{1.4cm}<{\centering}|m{0.8cm}<{\centering}}
    \hline
    \multicolumn{2}{c|}{Case} & Def & Atk & OC(\$) & \tabincell{c}{NP\\($MWh$)}\\
    \hline
    \multirow{2}*{1} &OGA&$PLine_{5}$ $PLine_{7}$&$PLine_{1}$ $PLine_{3}$&$2.764\times 10^5$&$56.1$\\
    \cline{2-6}
    &no OGA&$PLine_{3}$ $PLine_{4}$&$PLine_{1}$&$2.201\times 10^5$&$0$\\
    \hline
    \multirow{2}*{2} &OGA&$PLine_{1}$ $PLine_{2}$&$PLine_{4}$ $PLine_{5}$&$3.469\times 10^5$&$49.6$\\
    \cline{2-6}
    &no OGA&$PLine_{2}$ $PLine_{4}$&$PLine_{1}$&$2.969\times 10^5$&$0$\\
    \hline
    \multirow{2}*{3} &OGA&$PLine_{2}$ $PLine_{3}$&$PLine_{5}$ $PLine_{6}$&$4.405\times 10^5$&$432$\\
    \cline{2-6}
    &no OGA&$PLine_{3}$ $PLine_{4}$&$PLine_{1}$&$2.827\times 10^5$&$55.8$\\
    \hline
    \multirow{2}*{4} &OGA&$PLine_{2}$ $PLine_{6}$&$PLine_{3}$ $PLine_{4}$&$4.027\times 10^5$&$394$\\
    \cline{2-6}
    &no OGA&$PLine_{2}$ $PLine_{4}$&$PLine_{1}$&$2.194\times 10^5$&$7.25$\\
    \hline
  \end{tabular}
\end{table}

\subsection{Impact of $BigM_{obj}$}

\begin{figure}
\centering
  \includegraphics[width=0.35\textwidth]{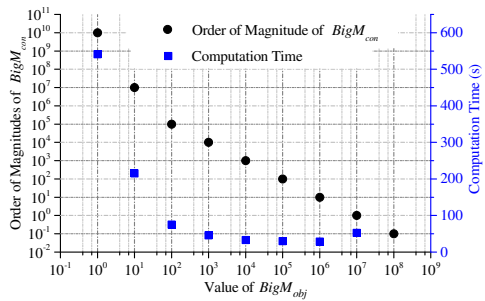}
  \caption{Value of $BigM$ and computation time.}
  \label{fig:BigM_obj}
\end{figure}

As mentioned in section III.C, the value of $BigM_{con}$ for unbounded dual variables will significantly impact the computational efficiency. Different from \cite{Arroyo_AD} and \cite{Yuan_DAD}, Fig. \ref{fig:BigM_obj} shows the active lower bound of $BigM_{con}$ as well as computational time along with the variation of $BigM_{obj}$. This shows, the active lower bound of the magnitude of $BigM_{con}$ keeps decreasing as $BigM_{obj}$ decreases, which confirms the analysis in section III.C. However, the computational time is non-monotonic with respect to $BigM_{obj}$. The reason may be that, at the very beginning, as $BigM_{obj}$ decreases, the range of $BigM_{con}$ shrinks, which would lead to a stronger MILP formulation and thus reduce the computational burden; however, as $BigM_{obj}$ keeps decreasing, the range of $BigM_{con}$ may become too small and step out of the numerical stability range, leading in turn to an ill-conditioned problem. Nevertheless, $BigM_{obj}$ has a relatively wide ideal range from the aspect of computational efficiency, which is [$10^4, 10^6]$ in this case, as Fig. \ref{fig:BigM_obj} indicates, demonstrating the promising prospect of practical application for the proposed method.

\subsection{39-Bus Power System with 20-Node Gas Network}

In this subsection, we apply the proposed model and algorithm to a larger test system, which comprises the IEEE 39-bus system and a modified version of the Belgian high-calorific 20-node gas network. It has 3 gas-fired units, 7 traditional units, 2 gas wells, 3 compressors, 4 gas storage tanks, 19 power loads and 9 gas loads. Similar to the settings in section IV. A, all the system branches are the targets of the attacker, including 46 power transmission lines, 21 passive pipelines, 3 active pipelines and 3 connection lines. Refer to \cite{Power6Gas7} for the topology, parameters, and unit commitment of the test system as well as the power and gas demand. We consider $T=2$ periods and $K=4$ for the Weymouth equation approximation and a 1$\%$ optimality gap as a trade-off between solution accuracy and computational burden. We select period 2 to 3, 9 to 10, 14 to 15 and 19 to 20, as four target time slots. Table \ref{Tab:Big_case_time} summarizes the average computation time under different combinations of defense and attack budget. The results show that the computation time increases rapidly as the defense and attack budget increase. However, when the defense and attack budget are limited, the computation time is acceptable, demonstrating the scalability of the proposed model and algorithm to moderately sized systems.

\begin{table}[h]
\footnotesize
  \centering
  \caption{Computation Time (s) under Different Budgets}
  \label{Tab:Big_case_time}
  \begin{tabular}{m{1.5cm}<{\centering}|m{0.8cm}<{\centering}|m{0.8cm}<{\centering}|m{0.8cm}<{\centering}|m{0.8cm}<{\centering}|m{0.8cm}<{\centering}}
    \hline
    \backslashbox{Def}{Atk} & 1 & 2 & 3 & 4 & 5\\
    \hline
    1 & 65 & 77 & 93 & 111 & 140     \\
    \hline
    2 & 155 & 172 & 312 & 551 & 673  \\
    \hline
    3 & 317 & 477 & 699 & 831 & 1145 \\
    \hline
    4 & 663 & 927 & 1401 & 2022 & 2712 \\
    \hline
    5 & 922 & 1477 & 1902 & 2638 & 3531 \\
    \hline
  \end{tabular}
\end{table}

\section{Conclusion}

The proposed methodology addresses the vulnerability of coupled gas-electric networks against malicious line interdictions. It provides the optimal strategies for preventive reinforcement, and increases the resilience of the energy supply, and decreases the operational cost of the interdependent energy systems. The mixed integer linear representation of the over-generation issue and the gas system constraints entails binary variables in the lower level problem, and prevents a traditional C\&CG algorithm being applied. A nested $C\&CG$ algorithm is adopted and the corresponding computational burden is reduced by choosing the proper value of $BigM_{con}$, which enhances the applicability of the proposed method. Simulation results corroborate the effectiveness of the proposed method, and suggest giving higher priority to protecting systems with radial topology.

\ifCLASSOPTIONcaptionsoff
  \newpage
\fi

\bibliographystyle{IEEEtran}
\bibliography{IEEEabrv,refs}

% biography section
%
% If you have an EPS/PDF photo (graphicx package needed) extra braces are
% needed around the contents of the optional argument to biography to prevent
% the LaTeX parser from getting confused when it sees the complicated
% \includegraphics command within an optional argument. (You could create
% your own custom macro containing the \includegraphics command to make things
% simpler here.)
%\begin{IEEEbiography}[{\includegraphics[width=1in,height=1.25in,clip,keepaspectratio]{mshell}}]{Michael Shell}
% or if you just want to reserve a space for a photo:

% You can push biographies down or up by placing
% a \vfill before or after them. The appropriate
% use of \vfill depends on what kind of text is
% on the last page and whether or not the columns
% are being equalized.

%\vfill

% Can be used to pull up biographies so that the bottom of the last one
% is flush with the other column.
%\enlargethispage{-5in}

% that's all folks
\end{document}